\documentstyle{article}
\author{
G.N.Parfionov\footnotemark[1], 
Yu.A.Romashev\footnotemark[1], 
R.R.Zapatrine\footnotemark[1]{\makebox[0.4em]{}}\footnotemark[2]
}
\date{} 

\title{Vectors and covectors in non-commutative setting}

\def\hom{{\rm Hom}} 
\def\om{\omega} 

\def\endproof{\hspace*{\fill}$\Box$}
\def\proof{\paragraph{Proof.}} 
\newtheorem{prop}{Proposition} 
 
\newtheorem{theo}{Theorem}

\begin{document}
\maketitle 
\addtocounter{footnote}{+1} 
\footnotetext{Department of Mathematics, SPb UEF, Griboyedova 30/32, 
191023, St-Petersburg, Russia (adderss for correspondence)} 
\addtocounter{footnote}{+1} 
\footnotetext{Division of Mathematics, 
Istituto per la Ricerca di Base, 
I-86075, Monteroduni (IS), Molise, Italy}
 
\begin{abstract} 
Following the guidelines of classical differential geometry the 
`building material' for the tensor calculus in non-commutative 
geometry is suggested. The algebraic account of moduli of vectors 
and covectors is carried out. 
\end{abstract} 

\section*{Introduction} 

The main feature of the mathematics of quantum mechanics is 
captured in its non-commutativity. Thus, in order to build a 
quantum theory of spacetime it would be reasonable to implement an 
amount of non-commutativity into the classical differential 
geometry and general relativity. To do it, we ought to use their 
algebraic formulation \cite{chev,geroch}. In particular, it was 
shown by R. Geroch \cite{geroch} that the entire content of general 
relativity can be reformulated in mere terms of the algebra 
${C^\infty(M)}$ 
of smooth functions on a spacetime manifold. 

The direct attempt to substitute the commutative algebra 
${C^\infty(M)}$ by 
a non-commutative one causes both mathematical and physical 
problems related with the ambiguity of the generalization of 
geometrical objects \cite{ps}. Besides that, there is a duality at 
the very starting point which can be briefly formulated as `what to 
begin with: covectors or vectors?' The first opportunity was 
investigated by M.Karoubi \cite{karoubi} and A.Connes 
\cite{connes}.  Rather, in our paper we shall deal with vectors as 
basic objects keeping ourselves closer to the conventional account 
of the differential geometry. It should pointed out that in the 
classical theory both accounts are equivalent while in the 
non-commutative environment this is no longer so and the resulting 
`non-commutative geometries' are different. 

We would like to outline the liaisons of our approach with the 
`French version' of non-commutative geometry 
\cite{connes,dubviol,karoubi}. For us, the starting object is a 
dual to the bimodule $\Omega^1$ rather than $\Omega^1$ itself. 
However, the consequence of our construction is the appearance of a 
kind of `ghosts', that is, the covectors which can not be expressed 
in terms of differential forms. 

\section{Differential algebras} \label{s1} 

Let ${\cal A}$ be an associative algebra, ${{\rm Der}({\cal A})}$ 
be the set of its 
derivatives, that is, the linear mappings $v:{\cal A}\to {\cal A}$ 
for which the Leibniz rule holds: 
\begin{equation}\label{f8} 
v(ab) = va \cdot b + a \cdot vb 
\end{equation} 

${{\rm Der}({\cal A})}$ is a vector space in the sense that any derivative 
multiplied by a number remains derivative. When ${\cal A}  = {C^\infty(M)}$ 
the space ${{\rm Der}({\cal A})}$ is formed by all 
smooth vector fields on $M$. 

Even in classical differential geometry ${{\rm Der}({\cal A})}$ 
considered vector space has infinite dimension. 
Enlarging the set of multiplicators from the numbers to the 
elements of ${\cal A}$ we make ${{\rm Der}({\cal A})}$ 
left ${\cal A}$-module which is always 
finitely generated. This module may not be free being however 
projective due to Swan theorem \cite{swan} (see the example below). 

Sometimes subspaces $V\subseteq {{\rm Der}({\cal A})}$ of derivatives 
rather than 
the whole ${{\rm Der}({\cal A})}$ are considered. In the classical 
setting it 
happens in the two following situations. The first one is when 
theories with a symmetry group are considered                   
and only the invariant vector fields are taken into account. The 
second one (gauge theories) is when fiber bundles are considered 
and the vector fields on the total space tangent to 
the fibers are specified. 

To pass to the non-commutative setting we begin with the 
{\sc differential algebra} being the couple $({\cal A},V)$ where 
${\cal A}$ is an associative algebra and 
$V\subseteq {{\rm Der}({\cal A})}$ is a linear 
subspace of ${{\rm Der}({\cal A})}$. As in the classical setting, 
to reduce the 
number of generators of $V$ we attempt to endow $V$ with the 
structure of a module. Let, for instance, $V = {{\rm Der}({\cal A})}$, 
$v\in {\cal A}$. 
As in the classical setting, defining for an arbitrary element 
$s\in {\cal A}$ the linear mapping $sv: {\cal A}\to {\cal A}$ as 
$(sv)a = s\cdot va$. 
Checking the Leibniz rule (\ref{f8}) for the mapping $sv$ 
\[ 
(sv)(ab) = sva\cdot b + a\cdot svb - [a,s]\cdot vb 
\] 
for non-commutative ${\cal A}$ we see that it may not hold in general. 
However for any $s\in {\cal Z}$ (the center of ${\cal A}$) the 
mapping $sv$ will be the element of ${{\rm Der}({\cal A})}$, 
therefore ${{\rm Der}({\cal A})}$ is always 
${\cal Z}$-module. That is why in the definition of the differential 
algebra we shall always require $V$ to be a ${\cal Z}$-submodule of 
${{\rm Der}({\cal A})}$. 

\section{The coupling procedure} \label{s2} 

A coupling procedure binding the vectors and covectors is the 
necessary condition to introduce such basic geometrical entities 
as, say, curvature. The notions of vectors and covectors are 
introduced in this section based on the coupling procedure borrowed 
from the classical differential geometry. 

Let $({\cal A},V)$ be a differential algebra. As in the classical 
setting, we shall call the elements of the module $V$ {\sc 
vectors}. Rather, the covectors are yet to be defined and there are 
{\em a priori} several ways to do it. Clearly, the coupling 
procedure needs at least an object to be coupled with. In classical 
differential geometry this object is uniquely defined being 
the dual ${\cal A}$-module to $V$. In the non-commutative case 
it can not be carried out since $V$ is not an ${\cal A}$-module 
(as it was shown above) and a thoroughful analysis of the 
available opportunities is needed. 

Recall that $V$ is nevertheless ${\cal Z}$-module and consider its 
standard algebraic dual 
\begin{equation}\label{f120} 
V^\ast = \hom_{{\cal Z}}(V,{\cal Z}) 
\end{equation} 
where $\hom_{{\cal Z}}$ means the set of all ${\cal Z}$-linear 
mappings. Thus the coupling is automatically defined for any $\om 
\in V^\ast, v \in V$  
\begin{equation}\label{f125} \label{coupling} 
<\om, v> = \om(v) 
\end{equation} 
It is by virtue of the definition of $<\om,v>$ that this form 
is ${\cal Z}$-linear by the second argument. Note that in the classical 
setting the introduced object is exactly the module of covectors.  
Thus it seems natural to call the elements of $V^{\dag}$ {\sc 
covectors} as well. 

However, in order to follow the geometrical guideline, it would be 
reasonable to require the differentials $da$ to be covectors, where 
as usually, 
\[ da(v) = va \] 
for any $a\in {\cal A}, v\in V$. This requirement is not in general 
compatible with the definition (\ref{f120}), since the values of 
$va$ may occur beyond the center ${\cal Z}$. 

To gather it we shall expand the range of the values of the form 
(\ref{f125}) to the whole ${\cal A}$. Thus the following dual object 
$V^{\dag}$ is suggested 
\begin{equation}\label{f145} 
V^{\dag} = \hom_{{\cal Z}}(V,{\cal A}) 
\end{equation} 
keeping the same definition of the coupling bracket (\ref{f125}). 

\begin{prop} The linear space $V^{\dag}$ is ${\cal A}$-bimodule with 
respect to the following action of ${\cal A}$: 
\begin{equation}\label{f160} 
(a\om b)(v) = a\cdot \om(v) \cdot b 
\end{equation} 
where $a,b\in {\cal A}, \om \in V^{\dag}, v\in V$
\end{prop} 

\proof It suffices to prove that $a \om b$ is ${\cal Z}$-linear \( (a\om 
b)(zv) = a\cdot z\om(v) \cdot b = z\cdot a\cdot \om(v) \cdot b \) 
since $z$ commutes with any element of ${\cal A}$. \endproof 

\paragraph{Corollary.} 1. This proposition enables the Leibniz rule 
to hold for all differentials:  \( d(ab) = da\cdot b + a\cdot db 
\). 

\noindent 2. The bilinear form $<\cdot,\cdot>$ (\ref{coupling}) is 
${\cal A}$-linear with respect to the first argument and 
${\cal Z}$-linear with 
respect to the second one. 

\begin{prop} The bilinear form (\ref{coupling}) is 
nondegenerate: 
\begin{equation}\label{f190} 
\begin{array}{rcccl} 
<\om,v>=0 &,& \forall v\in V &\Rightarrow& \om=0 \cr
<\om,v>=0 &,& \forall \om\in V^{\dag} &\Rightarrow& v=0 
\end{array} 
\end{equation} 
\end{prop} 

\paragraph{Sketch of proof.} The first implication holds by the 
definition of $V^{\dag}$. To prove the second it suffices to check it 
for all differentials $da$.  \endproof

\section{The second dual} \label{s3} 

The non-commutativity of the basic algebra ${\cal A}$ breakes the 
symmetry between the vectors and covectors: the vectors form a 
${\cal Z}$-module while the covectors are ${\cal A}$-bimodule. Moreover, the 
transition from $V$ to $V^{\dag}$ is not the conjugation of 
${\cal Z}$-moduli. It is asymmetry that compensates asymmetry. So, we 
take 
\begin{equation}\label{f210} 
V^{{\dag}{\dag}} = \hom_{{\cal A}}(V^{\dag},{\cal A}) 
\end{equation} 
(the set of all homomorphisms of ${\cal A}$-bimoduli) as the second dual 
to $V$. Now the symmetry is restored which is corroborated by the 
following proposition. 

\begin{prop} $V^{{\dag}{\dag}}$ is ${\cal Z}$-module with respect 
to the standard action of ${\cal Z}$: 
\[ (zw)(\om) = z\cdot w(\om) \] 
for $z\in {\cal Z}, w\in V^{{\dag}{\dag}}, \om \in V^{\dag}$. 
\end{prop} 

\paragraph{Proof} is obtained by direct checking the ${\cal A}$ 
linearity of $zw$ using the commutativity of the elements of ${\cal Z}$. 
\endproof 

In the general theory of moduli there is the canonical {\em 
homomorphism} from a module to its second dual. In our setting this 
homomorphism is even injective: 

\begin{prop} The canonical homomorphism $v\mapsto \hat{v}$ 
\begin{equation}\label{f230} 
\hat{v}(\om) = <\om,v> = \om(v) 
\end{equation} 
is the embedding $V\to V^{{\dag}{\dag}}$ . 
\end{prop} 

\paragraph{Proof} follows immediately from the nondegeneracy of the 
coupling form (\ref{f190}) \endproof 

Recall that in classical differential geometry 
$V \simeq V^{{\dag}{\dag}} = V^{\ast\ast}$ so the reflexivity 
always takes place. 

\section{Projectivity and reflexivity} \label{s4} 

To build the working tensor calculus including the trace (which is 
necessary, in particular, to form the Ricci tensor) we have to 
deal with reflexive moduli $V \simeq V^{{\dag}{\dag}}$. For general 
differential algebras $({\cal A},V)$ this may not hold. In this section 
we show that in the case when $V$ is a projective ${\cal Z}$-module (as 
it is always in classical differential geometry according to Swan's theorem 
\cite{swan}) the 
reflexivity of the module of vectors is guaranteed. 

\begin{theo} \label{th1}  Let $({\cal A},V)$ is a differential algebra 
such that the module $V$ is a projective finitely generated 
${\cal Z}$-module.  Then the canonical embedding (\ref{f230}) $V \to 
V^{{\dag}{\dag}}$ is isomorphism, that is, the module $V$ is reflexive.  
\end{theo} 

\proof We shall use the following definition of projectivity: there 
exist a set of generators $\{v_1,\ldots ,v_n\}$ in $V$, and a set 
of cogenerators $\{\om^1,\ldots ,\om^n\}$ in $V^\ast = 
\hom_{{\cal Z}}(V,{\cal Z})$ (sic!) such that for any $v\in V$ 
\begin{equation}\label{f265} 
v = \om^1(v)v_1 + \ldots +\om^n(v)v_n 
\end{equation} 
First prove that for any $\om\in V^{\dag}$ (rather than from $V^\ast$) 
\[ 
\om = \om^1\om(v_1) + \ldots + \om^n\om(v_n) 
\] 
which is obtained by calculation of the value of $\om$ on 
arbitrary $v\in V$ decomposed by (\ref{f265}). 

Now consider an arbitrary $w\in V^{{\dag}{\dag}}$ and find such $v\in V$ that 
$\hat{v}=w$. Introduce the following mapping $v:{\cal A}\to {\cal A}$: 
\[ 
v(a) = \sum w(\om^i)v_ia 
\] 
which is the element of $V$ since the coefficients $w(\om^i)$ are 
always in ${\cal Z}$: $\om^i(u)\in {\cal Z}$ for any $u\in V$, hence $\om^ia = 
a\om^i$ for all $a\in {\cal A}$, therefore $w(\om^i)a = aw(\om^i)$ since 
$w$ is homomorphism of ${\cal A}$-bimoduli, thus $v\in V$. 

Finally, calculating directly the value of $\hat{v}$ on arbitrary 
$\om \in V^{\dag}$ we obtain $\hat{v}(\om) = w(\om)$, which completes 
the proof. \endproof 

\section*{Concluding remarks} 

Some conditions for developing the tensor calculus in 
non-commutative setting were studied in this paper in order to make 
it applicable for the quantization of gravity. That is, we had to 
follow the guidelines provided by the Einstein's theory based on 
the classical differential geometry. In particular, we should take 
care of the correspondence principle so that our construction would 
really be a generalization of classical differential geometry. 
Passing to the non-commutative setting yielded us a sort of 
`ghosts' not existing in the classical theory:  it turns out that 
the module of vectors may not be reflexive, and the dual space 
contains something more than vectors. Although, it was shown 
(theorem \ref{th1}) that in the situations similar to classical 
ones these ghosts  do not exist and there are no interpretational 
problems in building the non-commutative version of differential 
geometry.


\begin{thebibliography}{99} 

\bibitem{chev} Chevalley, K., Th\'eorie des groupes de Lie,
Hermann, Paris, 1955

\bibitem{connes} Connes A., Noncommutative differential geometry, 
Hermann, Paris, 1989 

\bibitem{dubviol} Dubois-Violette, M., {\it D\'erivations et calcul
diff\'erentiel non-commutatif}, Comptes Rendus de l'Academie
des Sciences de Paris, ser. I,
{\bf 307}, 403, (1988)

\bibitem{geroch} Geroch, R. 
{\it Einstein Algebras}, 
Communications in Mathematical Physics, 
{\bf 26}, 
271, 
1972 

\bibitem{karoubi} Karoubi, M., {\it Homologie cyclique et
K-th\'eorie}, Ast\'erisque,
{\bf 149}, 1, (1987)

\bibitem{ps} G.N.Parfionov, R.R.Zapatrin, 
{\it Pointless Spaces in General Relativity}, 
International Journal of Theoretical Physics, 
{\bf 34}, 737, 
1995 (eprint gr-qc/9503048) 

\bibitem{swan} Swan, R.G., {\it Vector fields and projective
modules}, Transactions of the AMS,
{\bf 105},
264,
(1962)


\end{thebibliography}
\end{document}